\documentclass[10pt]{amsart} 

\usepackage{amssymb}

\newtheorem{thm1}{Theorem}[section]

\newtheorem{rem1}[thm1]{Remark}

\newtheorem{cor1}[thm1]{Corollary}

\newtheorem{prop1}[thm1]{Proposition}
\newtheorem{ex1}[thm1]{Example}
\newtheorem{not1}[thm1]{Notation}

\begin{document}

\title[]
{Projective closure of Gorenstein monomial curves and the Cohen-Macaulay property}
\author[]{Anargyros Katsabekis}
\address { Department of Mathematics, University of Ioannina, 45110 Ioannina, Greece} \email{katsampekis@uoi.gr}

\begin{abstract} Let $C({\bf a})$ be a Gorenstein non-complete intersection monomial curve in the 4-dimensional affine space. There is a vector ${\bf v} \in \mathbb{N}^{4}$ such that for every integer $m \geq 0$, the monomial curve $C({\bf a}+m{\bf v})$ is Gorenstein non-complete intersection whenever the entries of ${\bf a}+m{\bf v}$ are relatively prime. In this paper, we study the arithmetically Cohen-Macaulay property of the projective closure of $C({\bf a}+m{\bf v})$.
\end{abstract}
\maketitle

\section{Introduction}
Let ${\bf a}=(a_{1},\ldots,a_{n})$ be a sequence of positive integers with ${\rm gcd}(a_{1},\ldots,a_{n})=1$. Consider the polynomial ring $R=K[x_{1},\ldots,x_{n}]$ in $n$ variables over a field $K$. Let $C({\bf a})$ be the monomial curve in the $n$-dimensional affine space $\mathbb{A}^{n}(K)$ defined parametrically by $x_{1}=t^{a_1},\ldots,x_{n}=t^{a_n}$. The {\em toric ideal} of $C({\bf a})$, denoted by $I({\bf a})$, is the kernel of the $K$-algebra homomorphism $\phi: R \rightarrow K[t]$ given by $$\phi(x_{i})=t^{a_i} \ \ \textrm{for all} \ \ 1 \leq i \leq n.$$ It is known that $I({\bf a})$ is a prime ideal generated by binomials $x_{1}^{u_{1}} \cdots x_{n}^{u_{n}}-x_{1}^{v_{1}} \cdots x_{n}^{v_{n}}$ such that $u_{1}a_{1}+\cdots+u_{n}a_{n}=v_{1}a_{1}+\cdots+v_{n}a_{n}$, see \cite[Lemma 4.1]{Sturmfels95}. 

Let us assume that $a_{n}>a_{i}$ for all $i<n$. The {\em homogenization} of $I({\bf a})$ with respect to the variable $x_{0}$, denoted by $I^{h}({\bf a})$, is the kernel of the $K$-algebra homomorphism $\psi: R[x_{0}] \rightarrow K[s,t]$ given by $\psi(x_{0})=s^{a_n}$ and $\psi(x_{i})=t^{a_i}s^{a_{n}-a_{i}}$, for al $1 \leq i \leq n$. The projective monomial curve $\overline{C({\bf a})}$ in $\mathbb{P}^{n}$ is given parametrically by $[s:t] \mapsto [s^{a_n}:t^{a_1}s^{a_{n}-a_{1}}: \cdots: t^{a_{n-1}}s^{a_{n}-a_{n-1}}:t^{a_n}]$. If $K$ is algebraically closed, then $I^{h}({\bf a})$ is the vanishing ideal of $\overline{C({\bf a})}$. The projective monomial curve $\overline{C({\bf a})}$ is called {\em arithmetically Cohen-Macaulay} if its vanishing ideal $I^{h}({\bf a})$ is a Cohen-Macaulay ideal. Arithmetically Cohen–Macaulayness of projective monomial curves has been studied by many authors; see, for instance, \cite{CN}, \cite{HerSta}, \cite{RR1}, \cite{SSS}, \cite{SeSe}. 

Recently, there has been an increased interest in studying the behaviour of the Betti numbers of a shifted toric ideal $I({\bf a}+m{\bf u})$, where $m \geq 0$ is an integer and ${\bf u} \in \mathbb{N}^{n}$, see \cite{CGHN}, \cite{GS}, \cite{JS}, \cite{Katsa}, \cite{Vu}. That's mainly due to a conjecture of J. Herzog and H. Srinivasan saying that the Betti numbers of $I(a_{1}+m,\ldots,a_{n}+m)$ are eventually periodic in $m$. The conjecture was proved by T. Vu in \cite{Vu}. From Theorem 5.7 in \cite{Vu} it follows that the monomial curve $\overline{C(a_{1}+m,\ldots,a_{n}+m)}$ is arithmetically Cohen-Macaulay for all $m>>0$. 

Given a Gorenstein non-complete intersection monomial curve $C({\bf a})$ in $\mathbb{A}^{4}(K)$, P. Gimenez and H. Srinivasan \cite[Theorem 4]{GS} showed that there is a vector $(1,1,1,1) \neq {\bf v} \in \mathbb{N}^{4}$ such that for every $m \geq 0$, the monomial curve $C({\bf a}+m{\bf v})$ is Gorenstein non-complete intersection whenever the entries of ${\bf a}+m{\bf v}$ are relatively prime.

In this article our aim is to provide necessary and sufficient conditions for the arithmetically Cohen–Macaulayness of $\overline{C({\bf a}+m{\bf v})}$ by using a minimal generating set for $I({\bf a}+m{\bf v})$, see Theorems \ref{CM1} and \ref{CM6}. This information will allow us to check the arithmetically Cohen–Macaulay property by just computing a minimal generating set of the ideal. Our results apply also in the case $m=0$, thus providing necessary and sufficient conditions for the arithmetically Cohen–Macaulayness of $\overline{C({\bf a})}$. A basic step towards our goal is to find a Gr\"obner basis for the ideal $I({\bf a} +m{\bf v})$ with respect to a degree reverse lexicographic order, see Propositions \ref{Basic1}, \ref{Basic2}. We further determine a generating set for $I^{h}({\bf a}+m{\bf v})$, see Corollary \ref{BasicIndispensable1} and Corollary \ref{Indispensable4}. Finally we provide families of monomial curves $\overline{C({\bf a}+m{\bf v})}$ which are not arithmetically Cohen-Macaulay, see Examples \ref{Non-CM}, \ref{BasicExample}.

\section{Cohen-Macaulay criteria for $\overline{C({\bf a}+m{\bf v})}$}

In this section we first recall Bresinsky’s theorem, which gives the explicit description of $I({\bf a})$ when $C({\bf a})$ is a Gorenstein non-complete intersection monomial curve in $\mathbb{A}^{4}(K)$.

\begin{thm1} (\cite[Theorem 3]{Bresinsky75}) \label{Brebasic} Let $C({\bf a})$ be a monomial curve having the
parametrization $$x_1 = t^{a_1}, x_2 = t^{a_2}, x_3 = t^{a_3}, x_4 = t^{a_4}.$$
Then $C({\bf a})$ is a Gorenstein non-complete intersection monomial curve if and only if $I({\bf a})$ is minimally generated by the set
$$\{x_1^{d_1}- x_3^{d_{13}} x_4^{d_{14}}, x_{2}^{d_2}- x_{1}^{d_{21}}x_{3}^{d_{23}},  x_3^{d_{3}}-x_{2}^{d_{32}}x_{4}^{d_{34}}, x_{1}^{d_{41}}x_{2}^{d_{42}}-x_{4}^{d_4}, x_{2}^{d_{42}}x_3^{d_{13}}-x_{1}^{d_{21}}x_{4}^{d_{34}}\}$$ where the above binomials are unique up to isomorphism, $d_{ij}>0$ and also $$d_{1} =d_{21}+d_{41}, d_{2}= d_{32}+d_{42}, d_{3}=d_{13}+d_{23}, d_{4} =d_{14}+d_{34}.$$

\end{thm1}

\begin{rem1} \label{RemarkBasic} {\rm Bresinsky \cite[Theorem 4]{Bresinsky75} also showed that $C({\bf a})$ and $I({\bf a})$ are as in the previous theorem if and only if $a_{1}=d_{2}d_{4}d_{13}+d_{42}d_{14}d_{23}$, $a_{2}=d_{3}d_{4}d_{21}+d_{41}d_{34}d_{23}$, $a_{3}=d_{1}d_{2}d_{34}+d_{32}d_{21}d_{14}$, $a_{4}=d_{1}d_{3}d_{42}+d_{13}d_{32}d_{41}$ with ${\rm gcd}(a_{1},a_{2},a_{3},a_{4})=1$, $d_{i}>1$, $0<d_{ij}<d_{j}$ for $1 \leq i \leq 4$, and $d_{1} =d_{21}+d_{41}$, $d_{2}= d_{32}+d_{42}$, $d_{3}=d_{13}+d_{23}$, $d_{4} =d_{14}+d_{34}$.}
\end{rem1}

Let ${\bf a}=(a_{1},\ldots,a_{4})$ be a sequence of positive integers. A {\em binomial} $B=M-N \in I({\bf a})$ is called {\em  indispensable} of $I({\bf a})$ if every system of binomial generators of $I({\bf a})$ contains $B$ or $-B$, while a {\em monomial} $M$ is called {\em indispensable} of $I({\bf a})$ if every system of binomial generators of $I({\bf a})$ contains a binomial $B$ such that $M$ is a monomial of $B$. If $C({\bf a})$ is a Gorenstein non-complete intersection monomial curve, then $I({\bf a})$ is generated by the indispensable binomials, see \cite[Corollary 3.15]{KO}.

\begin{thm1} \label{GorBasic} (\cite[Theorem 4]{GS}, \cite[Theorem 8]{GS1}) Let $C({\bf a})$ be a Gorenstein non-complete intersection monomial curve in $\mathbb{A}^{4}(K)$. Consider the vector ${\bf v}=(v_{1},\ldots,v_{4}) \in \mathbb{N}^{4}$, where $v_{1}=d_{2}d_{3}-d_{23}d_{32}$, $v_{2}=d_{21}d_{3}+d_{23}d_{34}$, $v_{3}=d_{2}d_{34}+d_{21}d_{32}$ and $v_{4}=d_{2}d_{3}-d_{23}d_{32}$. For every $m \geq 0$ the monomial curve $C({\bf a}+m{\bf v})$ is Gorenstein non-complete intersection whenever the entries of ${\bf a}+m{\bf v}$ are relatively prime. Furthermore, $I({\bf a}+m{\bf v})$ is minimally generated by the binomials $f_{1}:=x_1^{d_1+m}- x_3^{d_{13}} x_4^{d_{14}+m}, f_{2}:=x_2^{d_2}- x_1^{d_{21}} x_3^{d_{23}}, f_{3}:=x_3^{d_3}- x_2^{d_{32}} x_4^{d_{34}}, f_{4}:= x_1^{d_{41}+m} x_2^{d_{42}}-x_4^{d_4+m}, f_{5}:=x_2^{d_{42}} x_3^{d_{13}}-x_1^{d_{21}}x_{4}^{d_{34}}$. 
\end{thm1}
Throughout this paper $C({\bf a})$ will always be a Gorenstein non-complete intersection monomial curve in $\mathbb{A}^{4}(K)$. Also we will use the notation introduced in Theorem \ref{GorBasic}.

\begin{rem1} {\rm If $d_{1} \geq d_{13}+d_{14}$ and $d_{21}+d_{34} \leq d_{42}+d_{13}$, then $d_{4} \leq d_{41}+d_{42}$ since $d_{14}+d_{34} \leq d_{1}-d_{13}+d_{34}=d_{21}+d_{41}-d_{13}+d_{34} \leq d_{41}+d_{42}$.}
\end{rem1}

\begin{prop1} \label{Basic1} Suppose that the following conditions hold.  \begin{enumerate}
\item $d_{1} \geq d_{13}+d_{14}$.
\item $d_{2} \geq d_{21}+d_{23}$.
\item $d_{3} \geq d_{32}+d_{34}$.
\item $d_{21}+d_{34} \leq d_{42}+d_{13}$.
\end{enumerate}
If the entries of ${\bf a}+m{\bf v}$ are relatively prime, then $G:=\{f_{i}| 1 \leq i \leq 5\}$ is the reduced Gr\"obner basis for $I({\bf a}+m{\bf v})$ with respect to the degree reverse lexicographic order $<$ with $x_{2}>x_{1}>x_{3}>x_{4}$.

\end{prop1}

\noindent \textbf{Proof.} We have that ${\rm in}_{<}(f_{1})=x_{1}^{d_{1}+m}$, ${\rm in}_{<}(f_{2})=x_{2}^{d_{2}}$, ${\rm in}_{<}(f_{3})=x_{3}^{d_{3}}$, ${\rm in}_{<}(f_{4})=x_{1}^{d_{41}+m}x_{2}^{d_{42}}$ and ${\rm in}_{<}(f_{5})=x_{2}^{d_{42}}x_{3}^{d_{13}}$. Since ${\rm in}_{<}(f_{1})$ and  ${\rm in}_{<}(f_{2})$ are relatively prime, we get $S(f_{1},f_{2}) \stackrel{G} {\longrightarrow} 0$. Similarly $S(f_{1},f_{3}) \stackrel{G} {\longrightarrow} 0$, $S(f_{1},f_{5}) \stackrel{G} {\longrightarrow} 0$, $S(f_{2},f_{3}) \stackrel{G} {\longrightarrow} 0$ and $S(f_{3},f_{4}) \stackrel{G} {\longrightarrow} 0$. We have that $S(f_{1},f_{4})=x_{1}^{d_{21}}x_{4}^{d_{4}+m}-x_{2}^{d_{42}}x_{3}^{d_{13}}x_{4}^{d_{14}+m} \stackrel{f_{5}} {\longrightarrow} 0$. Also $S(f_{2},f_{4})=x_{2}^{d_{32}}x_{4}^{d_{4}+m}-x_{1}^{d_{1}+m}x_{3}^{d_{23}} \stackrel{f_{1}} {\longrightarrow} x_{2}^{d_{32}}x_{4}^{d_{4}+m}-x_{3}^{d_3}x_{4}^{d_{14}+m} \stackrel{f_{3}} {\longrightarrow} 0$. Additionally $S(f_{2},f_{5})=x_{1}^{d_{21}}x_{3}^{d_3}-x_{1}^{d_{21}}x_{2}^{d_{32}}x_{4}^{d_{34}} \stackrel{f_{3}}{\longrightarrow} 0$. Furthermore, $S(f_{3},f_{5})=x_{1}^{d_{21}}x_{3}^{d_{23}}x_{4}^{d_{34}}-x_{2}^{d_{2}}x_{4}^{d_{34}} \stackrel{f_{2}}{\longrightarrow} 0$. Finally $S(f_{4},f_{5})=x_{1}^{d_{1}+m}x_{4}^{d_{34}}-x_{3}^{d_{13}}x_{4}^{d_{4}+m} \stackrel{f_{1}}{\longrightarrow} 0$. \hfill $\square$\\

In this paper we will make extensive use of the following result.

\begin{thm1} (\cite[Theorem 2.2]{HerSta}) \label{BasicStHer} Let ${\bf a}=(a_{1},\ldots,a_{4})$ be a sequence of positive integers with $a_{4}>a_{i}$ for all $i<4$. Denote $<$ any reverse lexicographic order on $R=K[x_{1},\ldots,x_{4}]$ such that $x_4$ is the smallest variable. The following conditions are equivalent: \begin{enumerate} \item the projective monomial curve $\overline{C({\bf a})}$ is arithmetically Cohen-Macaulay. \item $x_{4}$ does not divide any minimal generator of ${\rm in}_{<}(I({\bf a}))$.
\end{enumerate}
\end{thm1}

First we consider the case that $d_{2} \geq d_{21}+d_{23}$.

\begin{thm1} \label{CM1} Let $d_{2} \geq d_{21}+d_{23}$ and $a_{4}+mv_{4}>a_{i}+mv_{i}$ for all $i<4$. Suppose that the entries of ${\bf a}+m{\bf v}$ are relatively prime. Then $\overline{C({\bf a}+m{\bf v})}$ is arithmetically Cohen-Macaulay if and only if the following conditions hold: \begin{enumerate}
\item $d_{1} \geq d_{13}+d_{14}$;
\item $d_{3} \geq d_{32}+d_{34}$;
\item $d_{21}+d_{34} \leq d_{42}+d_{13}$.
\end{enumerate}

\end{thm1}

\noindent \textbf{Proof.} ($\Leftarrow$) Suppose that (1), (2) and (3) are true. By Proposition \ref{Basic1}, $G$ is the reduced Gr\"obner basis for $I({\bf a}+m{\bf v})$ with respect to the degree reverse lexicographic order $<$ with $x_{2}>x_{1}>x_{3}>x_{4}$. Since $x_4$ does not divide ${\rm in}_{<}(f_{i})$ for every $1 \leq i \leq 5$, we have, from Theorem \ref{BasicStHer}, that $\overline{C({\bf a}+m{\bf v})}$ is arithmetically Cohen-Macaulay.\\
 ($\Rightarrow$) Suppose that $\overline{C({\bf a}+m{\bf v})}$ is arithmetically Cohen-Macaulay. The binomials $f_i$, $1 \leq i \leq 5$, are indispensable of $I({\bf a}+m{\bf v})$, so they belong to the reduced Gr\"obner basis for $I({\bf a}+m{\bf v})$ with respect to $<$. But $x_4$ does not divide ${\rm in}_{<}(f_{i})$ for every $1 \leq i \leq 5$, thus ${\rm in}_{<}(f_{1})=x_{1}^{d_{1}+m}$, ${\rm in}_{<}(f_{3})=x_{3}^{d_{3}}$ and ${\rm in}_{<}(f_{5})=x_{2}^{d_{42}}x_{3}^{d_{13}}$. Therefore (1), (2) and (3) are true. \hfill $\square$\\

\begin{ex1} \label{Non-CM} {\rm Let ${\bf a}=(2a^{2}+a-2,2a+1,2a+3,2a^2+a-1)$, where $a \geq 2$ is an integer. Then $C({\bf a})$ is a Gorenstein non-complete intersection monomial curve and $I({\bf a})$ is minimally generated by $x_1^2-x_3^{a-1}x_4$, $x_2^{a+1}-x_1x_3$, $x_3^a-x_2x_4$, $x_1x_2^a-x_4^2$ and $x_2^ax_3^{a-1}-x_1x_4$. Here ${\bf v}=(a^{2}+a-1, a+1, a+2, a^{2}+a-1)$. Let $g$ be the greatest common divisor of the entries of ${\bf a}+m{\bf v}$. Then $g$ divides $(2a^2+a-1)-(2a^2+a-2)$, so $g=1$. By Theorem \ref{CM1} the monomial curve $\overline{C({\bf a}+m{\bf v})}$ is arithmetically Cohen-Macaulay if and only if $a=2$. Consequently the projective monomial curve $\overline{C(8+5m,5+3m,7+4m,9+5m)}$ is arithmetically Cohen-Macaulay for every $m \geq 0$. If $a>2$, then the monomial curve $\overline{C({\bf a}+m{\bf v})}$ is not arithmetically Cohen-Macaulay for every $m \geq 0$.}
\end{ex1}

Given a polynomial $f \in K[x_{1},\ldots,x_{4}]$, we shall denote by $f^{h}$ the homogenization of $f$ with respect to the variable $x_{0}$. The next corollary provides a generating set for $I^{h}({\bf a}+m{\bf v})$ under the conditions of Theorem \ref{CM1}.

\begin{cor1} \label{BasicIndispensable1} Let $d_{2} \geq d_{21}+d_{23}$ and $a_{4}+mv_{4}>a_{i}+mv_{i}$ for all $i<4$. Suppose that the entries of ${\bf a}+m{\bf v}$ are relatively prime. If $\overline{C({\bf a}+m{\bf v})}$ is arithmetically Cohen-Macaulay, then $G^{h}:=\{f_{i}^{h}| 1 \leq i \leq 5\}$ is the reduced Gr\"obner basis for $I^{h}({\bf a}+m{\bf v})$ with respect to the degree reverse lexicographic order $<_{0}$ on $R[x_{0}]$ extended from $R$ with $x_{2}>_{0}x_{1}>_{0} x_{3}>_{0}x_{4}>_{0}x_{0}$.
\end{cor1}

\noindent \textbf{Proof.}  By Proposition \ref{Basic1}, $G$ is the reduced Gr\"obner basis for $I({\bf a}+m{\bf v})$ with respect to $<$. From \cite[Lemma 1.1]{HerSta} the set $$G^{h}=\{x_1^{d_1+m}- x_{0}^{d_{1}-d_{13}-d_{14}}x_3^{d_{13}} x_4^{d_{14}+m}, x_2^{d_2}- x_{0}^{d_{2}-d_{21}-d_{23}}x_1^{d_{21}} x_3^{d_{23}},$$ $$x_3^{d_3}-x_{0}^{d_{3}-d_{32}-d_{34}} x_2^{d_{32}} x_4^{d_{34}}, x_1^{d_{41}+m} x_2^{d_{42}}-x_{0}^{d_{41}+d_{42}-d_{4}}x_4^{d_4+m},$$ $$x_2^{d_{42}} x_3^{d_{13}}-x_{0}^{d_{13}+d_{42}-d_{21}-d_{34}}x_1^{d_{21}}x_{4}^{d_{34}}\}$$ is the reduced Gr\"obner basis for $I^{h}({\bf a}+m{\bf v})$ with respect to $<_{0}$. \hfill $\square$\\

Next we consider the case that $d_{2}<d_{21}+d_{23}$.

\begin{not1} {\rm We shall denote by $w$ the minimum element of $$\{l \in \mathbb{N}|d_{1}+m-ld_{21} \leq 0\} \cup \{l \in \mathbb{N}|d_{3}-ld_{23} \leq 0\}.$$ Notice that $w \geq 2$.}
\end{not1}
Consider the binomials $$p_{i}=x_{2}^{(i+1)d_{2}-d_{32}}x_{3}^{d_{3}-(i+1)d_{23}}-x_{1}^{(i+1)d_{21}}x_{4}^{d_{34}}, \ 0 \leq i \leq w-2,$$ $$q_{i}=x_{1}^{d_{1}+m-(i+1)d_{21}}x_{2}^{(i+1)d_{2}-d_{32}}-x_{3}^{id_{23}}x_{4}^{d_{4}+m}, \ 0 \leq i \leq w-2,$$ and $$r= \begin{cases} x_{2}^{wd_{2}-d_{32}}-x_{1}^{wd_{21}}x_{3}^{wd_{23}-d_{3}}x_{4}^{d_{34}}, \ \ \ \  \ \ \ \ \  \ \ \ \  \textrm{if} \  d_{1}+m-wd_{21}> 0\\
x_{2}^{wd_{2}-d_{32}}-x_{1}^{wd_{21}-d_{1}-m}x_{3}^{(w-1)d_{23}}x_{4}^{d_{4}+m}, \ \textrm{if} \  d_{1}+m-wd_{21} \leq  0\\
 \end{cases}$$ Note that $p_{0}=f_{5}$ and $q_{0}=f_{4}$.

\begin{prop1} \label{Basic2} Let $d_{2}<d_{21}+d_{23}$ and suppose that the following conditions hold.
\begin{enumerate} \item $d_{1} \geq d_{13}+d_{14}$.
\item $d_{3} \geq d_{32}+d_{34}$.
\item $d_{21}+d_{34} \leq d_{13}+d_{42}$.
\item $(w-1)(d_{2}-d_{21}-d_{23})+d_{3}-d_{32}-d_{34} \geq 0$.
\item $(w-1)(d_{2}-d_{21}-d_{23})+d_{1}+d_{23}-d_{4}-d_{32} \geq 0$.
\item $w(d_{2}-d_{21}-d_{23})+d_{3}-d_{32}-d_{34} \geq 0$ when $d_{1}+m-wd_{21}> 0$.
\item $w(d_{2}-d_{21}-d_{23})+d_{1}+d_{23}-d_{4}-d_{32} \geq 0$ when $d_{1}+m-wd_{21} \leq 0$.

\end{enumerate}
If the entries of ${\bf a}+m{\bf v}$ are relatively prime, then $$T=G \cup \{p_{i}| 1 \leq i \leq w-2\} \cup \{q_{i}| 1 \leq i \leq w-2\} \cup \{r\}$$ is a Gr\"obner basis for $I({\bf a}+m{\bf v})$ with respect to the degree reverse lexicographic order $<$ on $R$ with $x_{2}>x_{1}>x_{3}>x_{4}$.
\end{prop1}

\noindent \textbf{Proof.} Notice that the polynomials $p_{i}$, $q_{i}$, $r$ are in the ideal $I({\bf a}+m{\bf v})$ by computation, for every $1 \leq i \leq w-2$. Here ${\rm in}_{<}(f_{1})=x_{1}^{d_{1}+m}$, ${\rm in}_{<}(f_{2})=x_{1}^{d_{21}}x_{3}^{d_{23}}$ and ${\rm in}_{<}(f_{3})=x_{3}^{d_{3}}$. For every $0 \leq i \leq w-2$ we have that $i+1 \leq w-1$, thus $(i+1)(d_{2}-d_{21}-d_{23})+d_{3}-d_{32}-d_{34} \geq (w-1)(d_{2}-d_{21}-d_{23})+d_{3}-d_{32}-d_{34}$ since $d_{2}-d_{21}-d_{23}<0$. So $(i+1)(d_{2}-d_{21}-d_{23})+d_{3}-d_{32}-d_{34} \geq 0$ and therefore $(i+1)d_{2}-d_{32}+d_{3}-(i+1)d_{23} \geq (i+1)d_{21}+d_{34}$. Thus ${\rm in}_{<}(p_{i})=x_{2}^{(i+1)d_{2}-d_{32}}x_{3}^{d_{3}-(i+1)d_{23}}$, for every $0 \leq i \leq w-2$. Since $i+1 \leq w-1$ for every $0 \leq i \leq w-2$, we have that $(i+1)(d_{2}-d_{21}-d_{23})+d_{1}+d_{23}-d_{4}-d_{32} \geq (w-1)(d_{2}-d_{21}-d_{23})+d_{1}+d_{23}-d_{4}-d_{32}$ and therefore $(i+1)(d_{2}-d_{21}-d_{23})+d_{1}+d_{23}-d_{4}-d_{32} \geq 0$. Thus $d_{1}-(i+1)d_{21}+(i+1)d_{2}-d_{32} \geq id_{23}+d_{4}$, so ${\rm in}_{<}(q_{i})=x_{1}^{d_{1}+m-(i+1)d_{21}}x_{2}^{(i+1)d_{2}-d_{32}}$ for every $0 \leq i \leq w-2$. Finally ${\rm in}_{<}(r)=x_{2}^{wd_{2}-d_{32}}$. Since ${\rm in}_{<}(f_{1})$ and  ${\rm in}_{<}(f_{3})$ are relatively prime, we deduce that $S(f_{1},f_{3}) \stackrel{T} {\longrightarrow} 0$. Similarly $S(f_{1},p_{i}) \stackrel{T} {\longrightarrow} 0$, for every $0 \leq i \leq w-2$, $S(f_{1},r) \stackrel{T} {\longrightarrow} 0$, $S(f_{2},r) \stackrel{T} {\longrightarrow} 0$, $S(f_{3},q_{i}) \stackrel{T} {\longrightarrow} 0$, for every $0 \leq i \leq w-2$, and $S(f_{3},r) \stackrel{T} {\longrightarrow} 0$. 

We have that $S(f_{1},f_{2})=x_{1}^{d_{41}+m}x_{2}^{d_{2}}-x_{3}^{d_{3}}x_{4}^{d_{14}+m} \stackrel{f_{4}} {\longrightarrow} x_{2}^{d_{32}}x_{4}^{d_{4}+m}-x_{3}^{d_{3}}x_{4}^{d_{14}+m}  \stackrel{f_{3}} {\longrightarrow} 0$. We will show that $S(f_{1}, q_{i}) \stackrel{T}{\longrightarrow} 0$, for every $0 \leq i \leq w-2$. It holds that $$S(f_{1},q_{i})=x_{1}^{(i+1)d_{21}}x_{3}^{id_{23}}x_{4}^{d_{4}+m}-x_{2}^{(i+1)d_{2}-d_{32}}x_{3}^{d_{13}}x_{4}^{d_{14}+m} \stackrel{f_{5}} {\longrightarrow}$$ $$ x_{1}^{(i+1)d_{21}}x_{3}^{id_{23}}x_{4}^{d_{4}+m}-x_{1}^{d_{21}}x_{2}^{id_{2}}x_{4}^{d_{4}+m} \stackrel{f_{2}} {\longrightarrow} 0.$$ 

Also $S(f_{2},f_{3})=x_{1}^{d_{21}}x_{2}^{d_{32}}x_{4}^{d_{34}}-x_{2}^{d_{2}}x_{3}^{d_{13}}  \stackrel{f_{5}} {\longrightarrow} 0$. Next we will show that $S(f_{2},p_{i}) \stackrel{T} {\longrightarrow} 0$, for every $0 \leq i \leq w-2$. If $0 \leq i \leq w-3$, then $S(f_{2},p_{i})=x_{1}^{(i+2)d_{21}}x_{4}^{d_{34}}-x_{2}^{(i+2)d_{2}-d_{32}}x_{3}^{d_{3}-(i+2)d_{23}} \stackrel{p_{i+1}} {\longrightarrow} 0$. Suppose that $i=w-2$ and let $d_{3}-wd_{23} \leq 0$. Then $S(f_{2},p_{w-2})=x_{1}^{wd_{21}}x_{3}^{wd_{23}-d_{3}}x_{4}^{d_{34}}-x_{2}^{wd_{2}-d_{32}}$. Suppose first that $d_{1}+m-wd_{21}>0$. Then $S(f_{2},p_{w-2})=-r$ and therefore $S(f_{2},p_{w-2}) \stackrel{r} {\longrightarrow} 0$. Suppose now that  $d_{1}+m-wd_{21} \leq 0$, then $S(f_{2},p_{w-2}) \stackrel{r} {\longrightarrow} x_{1}^{wd_{21}-d_{1}}x_{3}^{wd_{23}-d_{3}}x_{4}^{d_{34}} f_{1} \stackrel{f_{1}} {\longrightarrow} 0$. If $d_{3}-wd_{23}>0$, then $d_{1}+m-wd_{21} \leq 0$ by definition of $w$ and therefore $S(f_{2},p_{w-2})=x_{1}^{wd_{21}}x_{4}^{d_{34}}-x_{2}^{wd_{2}-d_{32}}x_{3}^{d_{3}-wd_{23}} \stackrel{r} {\longrightarrow} x_{1}^{wd_{21}}x_{4}^{d_{34}}-x_{1}^{wd_{21}-d_{1}-m}x_{3}^{d_{13}}x_{4}^{d_{4}+m}=x_{1}^{wd_{21}-d_{1}-m}x_{4}^{d_{34}} f_{1} \stackrel{f_{1}} {\longrightarrow} 0$. We will show that $S(f_{2},q_{i}) \stackrel{T} {\longrightarrow} 0$. If $0 \leq i \leq w-3$, then $S(f_{2},q_{i})=x_{1}^{d_{1}+m-(i+2)d_{21}}x_{2}^{(i+2)d_{2}-d_{32}}-x_{3}^{(i+1)d_{23}}x_{4}^{d_{4}+m} \stackrel{q_{i+1}} {\longrightarrow} 0$. Suppose that $i=w-2$.  If $d_{1}+m-wd_{21} \leq 0$, then $S(f_{2},q_{w-2})=x_{2}^{wd_{2}-d_{32}}-x_{1}^{wd_{21}-d_{1}-m}x_{3}^{(w-1)d_{23}}x_{4}^{d_{4}+m} \stackrel{r} {\longrightarrow} 0$. If $d_{1}+m-wd_{21}>0$, then $S(f_{2},q_{w-2})=x_{3}^{(w-1)d_{23}}x_{4}^{d_{4}+m}-x_{1}^{d_{1}+m-wd_{21}}x_{2}^{wd_{2}-d_{32}} \stackrel{r} {\longrightarrow} -x_{3}^{wd_{23}-d_{3}}x_{4}^{d_{34}}f_{1} \stackrel{f_{1}} {\longrightarrow} 0$. For every $0 \leq i \leq w-2$ we have that $$S(f_{3},p_{i})=x_{1}^{(i+1)d_{21}}x_{3}^{(i+1)d_{23}}x_{4}^{d_{34}}-x_{2}^{(i+1)d_{2}}x_{4}^{d_{34}} \stackrel{f_{2}} {\longrightarrow} 0.$$

Next we will show that $S(p_{i},p_{j}) \stackrel{T}{\longrightarrow} 0$, where $0 \leq i \leq w-2$ and $0 \leq j \leq w-2$ with $i<j$. We have that $S(p_{i},p_{j})=x_{1}^{(i+1)d_{21}}x_{4}^{d_{34}} (x_{1}^{(j-i)d_{21}}x_{3}^{(j-i)d_{23}}-x_{2}^{(j-i)d_{2}}) \stackrel{f_{2}}{\longrightarrow} 0$. We will show that $S(p_{i},r) \stackrel{T}{\longrightarrow} 0$. If $d_{1}+m-wd_{21}>0$, then $$S(p_{i},r)=x_{1}^{(i+1)d_{21}}x_{4}^{d_{34}}(x_{1}^{(w-i-1)d_{21}}x_{3}^{(w-i-1)d_{23}}-x_{2}^{(w-i-1)d_{2}}) \stackrel{f_{2}}{\longrightarrow} 0.$$ Suppose that $d_{1}+m-wd_{21} \leq 0$. Then $$S(p_{i},r)=x_{1}^{wd_{21}-d_{1}-m}x_{4}^{d_{34}}(x_{3}^{d_{3}+(w-i-2)d_{23}}x_{4}^{d_{14}+m}-x_{1}^{d_{1}+m-(w-i-1)d_{21}}x_{2}^{(w-i-1)d_{2}}) \stackrel{q_{w-i-2}}{\longrightarrow}$$ $$x_{1}^{wd_{21}-d_{1}-m}x_{3}^{(w-i-2)d_{23}}x_{4}^{d_{4}+m}(x_{3}^{d_3}-x_{2}^{d_{32}}x_{4}^{d_{34}}) \stackrel{f_{3}}{\longrightarrow} 0.$$

Next we will show that $S(q_{i},q_{j}) \stackrel{T}{\longrightarrow} 0$, where $0 \leq i \leq w-2$ and $0 \leq j \leq w-2$ with $i<j$. We have that $S(q_{i},q_{j})=x_{3}^{id_{23}}x_{4}^{d_{4}+m}(x_{1}^{(j-i)d_{21}}x_{3}^{(j-i)d_{23}}-x_{2}^{(j-i)d_{2}}) \stackrel{f_{2}}{\longrightarrow} 0$. We will prove that $S(q_{i},r) \stackrel{T}{\longrightarrow} 0$. If $d_{1}+m-wd_{21} \leq 0$, then $$S(q_{i},r)=x_{3}^{id_{23}}x_{4}^{d_{4}+m}(x_{1}^{(w-i-1)d_{21}}x_{3}^{(w-i-1)d_{23}}-x_{2}^{(w-i-1)d_{2}}) \stackrel{f_{2}}{\longrightarrow} 0.$$ Suppose that $d_{1}+m-wd_{21}> 0$. Then $$S(q_{i},r)=x_{1}^{d_{1}+m+(w-i-1)d_{21}}x_{3}^{wd_{23}-d_{3}}x_{4}^{d_{34}}-x_{2}^{(w-i-1)d_{2}}x_{3}^{id_{23}}x_{4}^{d_{4}+m} \stackrel{f_{1}}{\longrightarrow}$$ $$ x_{3}^{id_{23}}x_{4}^{d_{4}+m}(x_{1}^{(w-i-1)d_{21}}x_{3}^{(w-i-1)d_{23}}-x_{2}^{(w-i-1)d_{2}}) \stackrel{f_{2}}{\longrightarrow} 0.$$

Finally we show that $S(p_{i},q_{j}) \stackrel{T}{\longrightarrow} 0$, where $0 \leq i \leq w-2$ and $0 \leq j \leq w-2$. Suppose that $i<j$. Then $$S(p_{i},q_{j})=x_{3}^{d_{3}+(j-i-1)d_{23}}x_{4}^{d_{4}+m}-x_{1}^{d_{1}+m-(j-i)d_{21}}x_{2}^{(j-i)d_{2}}x_{4}^{d_{34}} \stackrel{f_{3}}{\longrightarrow}$$ $$x_{2}^{d_{32}}x_{4}^{d_{34}}(x_{3}^{(j-i-1)d_{23}}x_{4}^{d_{4}+m}-x_{1}^{d_{1}+m-(j-i)d_{21}}x_{2}^{(j-i)d_{2}-d_{32}}) \stackrel{q_{j-i-1}}{\longrightarrow} 0.$$ Assume that $i \geq j$, then $$S(p_{i},q_{j})=x_{2}^{(i-j)d_{2}}x_{3}^{d_{3}-(i-j+1)d_{23}}x_{4}^{d_{4}+m}-x_{1}^{d_{1}+m+(i-j)d_{21}}x_{4}^{d_{34}} \stackrel{f_{1}}{\longrightarrow}$$ $$x_{3}^{d_{3}-(i-j+1)d_{23}}x_{4}^{d_{4}+m}(x_{2}^{(i-j)d_{2}}-x_{1}^{(i-j)d_{21}}x_{3}^{(i-j)d_{23}}) \stackrel{f_{2}}{\longrightarrow} 0. \ \ \ \ \ \square $$

\begin{thm1} \label{CM6} Let $d_{2}<d_{21}+d_{23}$. Suppose that $a_{4}+mv_{4}>a_{i}+mv_{i}$ for all $i<4$ and the entries of ${\bf a}+m{\bf v}$ are relatively prime. Then $\overline{C({\bf a}+m{\bf v})}$ is arithmetically Cohen-Macaulay if and only if the following conditions hold: \begin{enumerate} \item $d_{1} \geq d_{13}+d_{14}$;
\item $d_{3} \geq d_{32}+d_{34}$;
\item $d_{21}+d_{34} \leq d_{13}+d_{42}$;
\item $(w-1)(d_{2}-d_{21}-d_{23})+d_{3}-d_{32}-d_{34} \geq 0$;
\item $(w-1)(d_{2}-d_{21}-d_{23})+d_{1}+d_{23}-d_{4}-d_{32} \geq 0$;
\item $w(d_{2}-d_{21}-d_{23})+d_{3}-d_{32}-d_{34} \geq 0$ when $d_{1}+m-wd_{21}> 0$;
\item $w(d_{2}-d_{21}-d_{23})+d_{1}+d_{23}-d_{4}-d_{32} \geq 0$ when $d_{1}+m-wd_{21} \leq 0$.
\end{enumerate}

\end{thm1}

\noindent \textbf{Proof.} ($\Leftarrow$) Suppose that (1), (2), (3), (4), (5), (6) and (7) are true. By Proposition \ref{Basic2}, $T$ is a Gr\"obner basis for $I({\bf a}+m{\bf v})$ with respect to the degree reverse lexicographic order $<$ with $x_{2}>x_{1}>x_{3}>x_{4}$. Since $x_4$ does not divide the initial monomial of any binomial in $T$, we have, from Theorem \ref{BasicStHer}, that $\overline{C({\bf a}+m{\bf v})}$ is arithmetically Cohen-Macaulay.\\ 
($\Rightarrow$) Suppose that $\overline{C({\bf a}+m{\bf v})}$ is arithmetically Cohen-Macaulay. Since the binomials $f_{1}$, $f_{3}$ and $f_{5}$ are indispensable of $I({\bf a}+m{\bf v})$, they belong to the reduced Gr\"obner basis of $I({\bf a}+m{\bf v})$ with respect to $<$. Consequently (1), (2) and (3) are true. Suppose that $(w-1)(d_{2}-d_{21}-d_{23})+d_{3}-d_{32}-d_{34}<0$. Then the binomial $p_{w-2}$ belongs to $I({\bf a}+m{\bf v})$ and ${\rm in}_{<}(p_{w-2})=x_{1}^{(w-1)d_{21}}x_{4}^{d_{34}}$. Note that $(w-1)d_{21}<d_{1}+m$ by definition of $w$. Since $\overline{C({\bf a}+m{\bf v})}$ is arithmetically Cohen-Macaulay, there is a binomial $B \in I({\bf a}+m{\bf v})$ such that ${\rm in}_{<}(B)$ divides $x_{1}^{(w-1)d_{21}}$, a contradiction to the fact that the monomial $x_{1}^{d_{1}+m}$ is indispensable. Assume that $(w-1)(d_{2}-d_{21}-d_{23})+d_{1}+d_{23}-d_{4}-d_{32}<0$. Then the binomial $q_{w-2}$ belongs to $I({\bf a}+m{\bf v})$ and ${\rm in}_{<}(q_{w-2})=x_{3}^{(w-2)d_{23}}x_{4}^{d_{4}+m}$. Note that $(w-2)d_{23}<d_{3}$ by definition of $w$. Since $\overline{C({\bf a}+m{\bf v})}$ is arithmetically Cohen-Macaulay, there is a binomial $u \in I({\bf a}+m{\bf v})$ such that ${\rm in}_{<}(u)$ divides $x_{3}^{(w-2)d_{23}}$, a contradiction to the fact that the monomial $x_{3}^{d_{3}}$ is indispensable. Suppose that $d_{1}+m-wd_{21}> 0$ and let $w(d_{2}-d_{21}-d_{23})+d_{3}-d_{32}-d_{34}<0$. Then the binomial $r$  belongs to $I({\bf a}+m{\bf v})$ and ${\rm in}_{<}(r)=x_{1}^{wd_{21}}x_{3}^{wd_{23}-d_{3}}x_{4}^{d_{34}}$. Since $\overline{C({\bf a}+m{\bf v})}$ is arithmetically Cohen-Macaulay, there is a binomial $B \in I({\bf a}+m{\bf v})$ such that ${\rm in}_{<}(B)$ divides $x_{1}^{wd_{21}}x_{3}^{wd_{23}-d_{3}}$. Thus $B=x_{1}^{b_1}x_{3}^{b_3}-x_{1}^{c_1}x_{2}^{b_2}x_{3}^{c_3}x_{4}^{b_4}$ with ${\rm in}_{<}(B)=x_{1}^{b_1}x_{3}^{b_3}$, where $b_{2} \neq 0$ or/and $b_{4} \neq 0$, $b_{1} \leq wd_{21}$ and $b_{3} \leq wd_{23}-d_{3}$. Note that $wd_{21}<d_{1}+m$ and $wd_{23}-d_{3}<d_{23}$ by definition of $w$, so $b_{1}<d_{1}+m$ and $b_{3}<d_{23}$. Moreover $b_{3}<d_{3}$ since $d_{23}<d_{3}$. But $B \in I({\bf a}+m{\bf v})$ and $I({\bf a}+m{\bf v})$ is generated by the set $\{f_{1},\ldots,f_{5}\}$, therefore there are polynomials $A_{i} \in K[x_{1},\ldots,x_{4}]$, $1 \leq i \leq 5$, such that $$B=A_{1}(x_1^{d_1+m}-x_3^{d_{13}} x_4^{d_{14}+m})+A_{2}(x_2^{d_2}- x_1^{d_{21}} x_3^{d_{23}})+A_{3}(x_3^{d_3}- x_2^{d_{32}} x_4^{d_{34}})+$$ $$A_{4}(x_1^{d_{41}+m} x_2^{d_{42}}-x_4^{d_4+m})+A_{5}(x_2^{d_{42}} x_3^{d_{13}}-x_1^{d_{21}}x_{4}^{d_{34}}).$$ In the previous equation let $x_{2}=0$ and $x_{4}=0$. Then $x_{1}^{b_1}x_{3}^{b_3}=A_{1}(x_{1},0,x_{3},0)x_1^{d_1+m}-A_{2}(x_{1},0,x_{3},0)x_1^{d_{21}} x_3^{d_{23}}+A_{3}(x_{1},0,x_{3},0)x_3^{d_3}$, thus ${\rm in}_{<}(B)$ is divided by at least one of the monomials $x_{1}^{d_{1}+m}$, $x_{1}^{d_{21}}x_{3}^{d_{23}}$ and $x_{3}^{d_{3}}$, a contradiction. Suppose that $d_{1}+m-wd_{21} \leq 0$ and let $w(d_{2}-d_{21}-d_{23})+d_{1}+d_{23}-d_{4}-d_{32}< 0$. Then $r$  belongs to $I({\bf a}+m{\bf v})$ and ${\rm in}_{<}(r)=x_{1}^{wd_{21}-d_{1}-m}x_{3}^{(w-1)d_{23}}x_{4}^{d_{4}+m}$. Since $\overline{C({\bf a}+m{\bf v})}$ is arithmetically Cohen-Macaulay, there is a binomial $u \in I({\bf a}+m{\bf v})$ such that ${\rm in}_{<}(u)$ divides $x_{1}^{wd_{21}-d_{1}-m}x_{3}^{(w-1)d_{23}}$. Thus ${\rm in}_{<}(u)=x_{1}^{e_1}x_{3}^{e_3}$ where $e_{1} \leq wd_{21}-d_{1}-m$ and $e_{3} \leq (w-1)d_{23}$. Note that $wd_{21}-d_{1}-m<d_{21}$ and $(w-1)d_{23}<d_{3}$ by definition of $w$, so $e_{1}<d_{21}$ and $e_{3}<d_{3}$. Moreover $e_{1}<d_{1}+m$ since $d_{21}<d_{1}+m$. But $u \in I({\bf a}+m{\bf v})$ and $I({\bf a}+m{\bf v})$ is generated by the set $\{f_{1},\ldots,f_{5}\}$, therefore ${\rm in}_{<}(u)$ is divided by at least one of the monomials $x_{1}^{d_{1}+m}$, $x_{1}^{d_{21}}x_{3}^{d_{23}}$ and $x_{3}^{d_{3}}$, a contradiction. \hfill $\square$\\

\begin{ex1} \label{BasicExample}{\rm Let ${\bf a}=(19,29,26,43)$, then $C({\bf a})$ is a Gorenstein non-complete intersection monomial curve and $I({\bf a})$ is minimally generated by $x_{1}^{5}-x_{3}^{2}x_{4}, x_{2}^{4}-x_{1}^{2}x_{3}^{3}, x_{3}^{5}-x_{2}^{3}x_{4}, x_{1}^{3}x_{2}-x_{4}^{2}, x_{2}x_{3}^{2}-x_{1}^{2}x_{4}$. Here ${\bf v}=(11,13,10,11)$, $w=2$ and $d_{1}+m-wd_{21}=m+1>0$. Suppose first that $m \leq 7$, then $43+11m={\rm max}\{19+11m,29+13m,26+10m,43+11m\}$. We have that $w(d_{2}-d_{21}-d_{23})+d_{3}-d_{32}-d_{34}=-1<0$, so from Theorem \ref{CM6} the monomial curve $\overline{C({\bf a}+m{\bf v})}$ is not arithmetically Cohen-Macaulay for $m \leq 7$ whenever ${\rm gcd}(19+11m,29+13m,26+10m,43+11m)=1$. Assume that $m>7$, then $29+13m={\rm max}\{19+11m,29+13m,26+10m,43+11m\}$. Let ${\bf b}=(b_{1},\ldots,b_{4})$ with $b_{1}=26+10m$, $b_{2}=19+11m$, $b_{3}=43+11m$, $b_{4}=29+13m$. Note that $b_{4}={\rm max}\{b_{1},\ldots,b_{4}\}$. Then $I({\bf b})$ is minimally generated by $x_{1}^{5}-x_{3}x_{4}^{3}, x_{2}^{5+m}-x_{1}^{2}x_{3}^{1+m}, x_{3}^{2+m}-x_{2}^{3+m}x_{4}, x_{1}^{3}x_{2}^{2}-x_{4}^{4}, x_{2}^2x_{3}-x_{1}^{2}x_{4}$. Since $x_{3}^{2+m}-x_{2}^{3+m}x_{4} \in I({\bf b})$, we have, from Theorem \ref{CM1}, that for every $m>7$ the monomial curve $\overline{C({\bf b})}$ is not arithmetically Cohen-Macaulay whenever ${\rm gcd}(b_{1}, b_{2}, b_{3}, b_{4})=1$.}
\end{ex1}

\begin{ex1} \label{BasicExample11}{\rm Let ${\bf a}=(a^{2}+5a,7a+6,6a+1,3a^{2}+3a-2)$ where $a>2$ is an even integer. Then $C({\bf a})$ is a Gorenstein non-complete intersection monomial curve and $I({\bf a})$ is minimally generated by $x_{1}^{3}-x_{3}^{2}x_{4}, x_{2}^{a}-x_{1}x_{3}^{a}, x_{3}^{a+2}-x_{2}x_{4}^2, x_{1}^{2}x_{2}^{a-1}-x_{4}^{3}, x_{2}^{a-1}x_{3}^{2}-x_{1}x_{4}^{2}$. Here ${\bf v}=(a^{2}+a,3a+2,2a+1,a^{2}+a)$, $w=2$ and $d_{1}+m-wd_{21}=m+1 >0$. It holds that $w(d_{2}-d_{21}-d_{23})+d_{3}-d_{32}-d_{34}=a-3>0$ and $(w-1)(d_{2}-d_{21}-d_{23})+d_{1}+d_{23}-d_{4}-d_{32}=a-2>0$, so, from Theorem \ref{CM6}, for every $m \geq 0$ the monomial curve $\overline{C({\bf a}+m{\bf v})}$ is arithmetically Cohen-Macaulay whenever the entries of ${\bf a}+m{\bf v}$ are relatively prime.}
\end{ex1}

\begin{ex1} {\rm Let ${\bf a}=(1191, 1239, 582, 2303)$, then $C({\bf a})$ is a Gorenstein non-complete intersection monomial curve and $I({\bf a})$ is minimally generated by $x_1^{16}-x_3^{9}x_4^6$, $x_2^{11}-x_1^{9}x_3^5$, $x_3^{14}-x_2x_4^{3}$, $x_1^{7}x_2^{10}-x_4^9$ and $x_2^{10}x_3^{9}-x_1^9x_4^{3}$. Notice that ${\bf v}=(149,141,42,149)$. Suppose first that $m \leq 2$, then $w=2$ and $d_{1}+m-wd_{21}=m-2 \leq 0$. We have that $w(d_{2}-d_{21}-d_{23})+d_{1}+d_{23}-d_{4}-d_{32}=5>0$ and $(w-1)(d_{2}-d_{21}-d_{23})+d_{3}-d_{32}-d_{34}=7>0$. By Theorem \ref{CM6}, for every $m \leq 2$ the monomial curve $\overline{C({\bf a}+m{\bf v})}$ is arithmetically Cohen-Macaulay whenever the entries of ${\bf a}+m{\bf v}$ are relatively prime. Suppose now that $m>2$, then $w=3$ and $d_{1}+m-wd_{21}=m-11$. Assume that $m \leq 11$, then $d_{1}+m-wd_{21} \leq 0$. We have that $w(d_{2}-d_{21}-d_{23})+d_{1}+d_{23}-d_{4}-d_{32}=2>0$ and $(w-1)(d_{2}-d_{21}-d_{23})+d_{3}-d_{32}-d_{34}=4>0$. By Theorem \ref{CM6}, the monomial curve $\overline{C({\bf a}+m{\bf v})}$ is arithmetically Cohen-Macaulay for every $2<m \leq 11$ whenever the entries of ${\bf a}+m{\bf v}$ are relatively prime. Suppose that $m>11$, then $d_{1}+m-wd_{21}>0$. We have that $w(d_{2}-d_{21}-d_{23})+d_{3}-d_{32}-d_{34}=1>0$  and $(w-1)(d_{2}-d_{21}-d_{23})+d_{1}+d_{23}-d_{4}-d_{32}=5>0$. By Theorem \ref{CM6}, the monomial curve $\overline{C({\bf a}+m{\bf v})}$ is arithmetically Cohen-Macaulay for every $m>11$ whenever the entries of ${\bf a}+m{\bf v}$ are relatively prime.}
\end{ex1}

Finally we provide a generating set for $I^{h}({\bf a}+m{\bf v})$ under the conditions of Theorem \ref{CM6}.

\begin{cor1} \label{Indispensable4} Let $d_{2}<d_{21}+d_{23}$. Suppose that $a_{4}+mv_{4}>a_{i}+mv_{i}$ for all $i<4$ and the entries of ${\bf a}+m{\bf v}$ are relatively prime. If $\overline{C({\bf a}+m{\bf v})}$ is arithmetically Cohen-Macaulay, then $T^{h}=G^{h} \cup \{p_{i}^{h}| 1 \leq i \leq w-2\} \cup \{q_{i}^{h}| 1 \leq i \leq w-2\} \cup \{r^{h}\}$  is a Gr\"obner basis for $I^{h}({\bf a}+m{\bf v})$ with respect to the degree reverse lexicographic order $<_{0}$ on $R[x_{0}]$ extended from $R$ with $x_{2}>_{0}x_{1}>_{0} x_{3}>_{0}x_{4}>_{0}x_{0}$.

\end{cor1}

\noindent \textbf{Proof.} We have that $$r^{h}= x_{2}^{wd_{2}-d_{32}}-x_{0}^{w(d_{2}-d_{21}-d_{23})+d_{1}+d_{23}-d_{4}-d_{32}}x_{1}^{wd_{21}-d_{1}-m}x_{3}^{(w-1)d_{23}}x_{4}^{d_{4}+m}$$ when $d_{1}+m-wd_{21} \leq  0$, while 
$$r^{h}=x_{2}^{wd_{2}-d_{32}}-x_{0}^{w(d_{2}-d_{21}-d_{23})+d_{3}-d_{32}-d_{34}}x_{1}^{wd_{21}}x_{3}^{wd_{23}-d_{3}}x_{4}^{d_{34}}$$ when $d_{1}+m-wd_{21}> 0$. From Proposition \ref{Basic2} the set $T$ is a Gr\"obner basis for $I({\bf a}+m{\bf v})$ with respect to $<$. By \cite[Lemma 1.1]{HerSta}, $T^{h}=\{f_{i}^{h}| 1 \leq i \leq 5\} \cup \{r^{h}\} \cup \{x_{2}^{(i+1)d_{2}-d_{32}}x_{3}^{d_{3}-(i+1)d_{23}}-x_{0}^{(i+1)(d_{2}-d_{21}-d_{23})+d_{3}-d_{32}-d_{34}}x_{1}^{(i+1)d_{21}}x_{4}^{d_{34}} | 1 \leq i \leq w-2\}  \cup \{x_{1}^{d_{1}+m-(i+1)d_{21}}x_{2}^{(i+1)d_{2}-d_{32}}-x_{0}^{(i+1)(d_{2}-d_{21}-d_{23})+d_{1}+d_{23}-d_{4}-d_{32}}x_{3}^{id_{23}}x_{4}^{d_{4}+m} | 1 \leq i \leq w-2\}$ is a Gr\"obner basis for $I^{h}({\bf a}+m{\bf v})$ with respect to  $<_{0}$. \hfill $\square$\\

\textbf{Acknowledgment.} The author would like to thank the referee for careful reading and helpful comments.


\begin{thebibliography}{50}


	
\bibitem{Bresinsky75} H. Bresinsky. {\em Symmetric semigroups of integers generated by 4 elements}, Manuscripta Math. {\bf 17} (3) (1975), 205–219. 
 \bibitem{CN} M. P. Cavaliere, G. Niesi. {\em On monomial curves and Cohen-Macaulay type}, Manuscripta Math. {\bf 42} (1983), no. 2-3, 147–159.
\bibitem{CGHN} R. Conaway, F. Gotti, J. Horton, C. O’Neill, R. Pelayo, M. Pracht, B. Wissman. {\em Minimal presentations of shifted numerical monoids}, Internat. J. Algebra Comput. {\bf 28} (2018), 53-68.
\bibitem{GS} P. Gimenez, H. Srinivasan. {\em A note on Gorenstein monomial curves},  Bull. Braz. Math. Soc. (N.S.) {\bf 45} (2014), no. 4, 671–678.
\bibitem{GS1} P. Gimenez, H. Srinivasan. {\em Structure of some numerical semigroup rings}, Banach Center Publications {\bf 121} (2020), 53-61.
\bibitem{HerSta} J. Herzog, D. Stamate. {\em Cohen–Macaulay criteria for projective monomial curves via Gr\"obner bases}, Acta Math. Vietnam. {\bf 44} (2019), 51–64.
\bibitem{JS} A.V. Jayanthan, H. Srinivasan. {\em Periodic occurence of complete intersection monomial curves}, Proc. Amer. Math. Soc. {\bf 141} (2013), 4199–4208.
\bibitem{Katsa}  A. Katsabekis. {\em Complete intersection monomial curves and the Cohen-Macaulayness of their tangent cones}, Algebra Colloq. {\bf 26} (2019), no. 4, 629–642.
\bibitem{KO} A. Katsabekis, I. Ojeda. {\em An indispensable classification of monomial curves in $\mathbb{A}^{4}(\mathbb{K})$}, Pacific J. Math. {\bf 268} (2014), 95–116. 

\bibitem{RR1} L. Reid, L. G. Roberts. {\em Non-Cohen-Macaulay projective monomial curves}, J. Algebra {\bf 291} (2005), no. 1, 171–186. 

\bibitem{SSS} J. Saha, I. Sengupta, P. Srivastava. {\em Betti sequence of the projective closure of affine monomial curves},  J. Symbolic Comput. {\bf 119} (2023), 101–111.
\bibitem{SeSe} T. Se, S. J. Grant. {\em The Cohen-Macaulay property of affine semigroup rings in dimension 2},  Comm. Algebra {\bf 47} (2019), no. 7, 2979–2994. 

\bibitem{Vu} T. Vu. {\em Periodicity of Betti numbers of monomial curves}, J. Algebra {\bf 418} (2014), 66–90.


\bibitem{Sturmfels95} B. Sturmfels. {\em Gr\"obner Bases and Convex Polytopes}, University Lecture Series Vol. 8. Providence, RI, USA: American Mathematical Society, 1996.
	
\end{thebibliography}
\end{document}